\documentclass[a4paper, 12pt]{article}
\usepackage{amssymb}

\usepackage{amsmath,amsthm}

\begin{document}

\title{Flocking Control of Groups of Mobile Autonomous Agents via Local Feedback  \thanks{Supported by the National Natural Science Foundation of China (No. 10372002 and No.
60274001) and the National Key Basic Research and Development
Program (No. 2002CB312200).}}
\author{Long Wang$^1$, Hong Shi$^1$, Tianguang Chu$^2$, Weicun Zhang$^3$, Lin
Zhang$^2$\\\\
$^1$Intelligent Control Laboratory, Center for Systems and Control,\\
Department of Mechanics and Engineering Science,\\
Peking University, Beijing 100871, P. R. China\\\\
$^2$Computer Science Department, Naval Postgraduate School,\\
Monterey, CA93943, USA\\\\
$^3$Automation Department, University of Science and Technology Beijing,\\
Beijing 100083, P. R. China}
\date{}
\maketitle


\textbf{Abstract}-- This paper considers a group of mobile
autonomous agents moving in Euclidean space with point mass
dynamics. We introduce a set of coordination control laws that
enable the group to generate the desired stable flocking motion.
The control laws are a combination of attractive/repulsive and
alignment forces. By using the control laws, all agent velocities
asymptotically approach the desired velocity, collisions can be
avoided between agents, and the final tight formation minimizes
all agent potentials. Moreover, we prove that the velocity of the
center of mass is always equal to the desired velocity or
exponentially converges to the desired value. Furthermore, we
study the motion of the group when the velocity damping is taken
into account. In this case, we can properly modify the control
laws to generate the same stable flocking motion. Finally, for the
case that not all agents know the desired final velocity, we show
that the desired flocking motion can still be obtained. Numerical
simulations are worked out to illustrate our theoretical results.


\section{Introduction}

Flocking motion can be found everywhere in nature, e.g., flocking
of birds, schooling of fish, and swarming of bees. Such collective
behavior has certain advantages such as avoiding predators and
increasing the chance of finding food. The study of collective
emergent behavior of multiple mobile autonomous agents has
attracted much attention in many fields such as biology, physics,
robotics and control engineering. Understanding the mechanisms and
operational principles in them can provide useful ideas for
developing distributed cooperative control and coordination of
multiple mobile autonomous agents/robots. Recently, distributed
control/coordination of the motion of multiple dynamic
agents/robots has emerged as a topic of major interest \cite{N. E.
Leonard and E. Fiorelli}--\cite{John H. Reif; H. Wang}. This is
partly due to recent technological advances in communication and
computation, and wide applications of multi-agent systems in many
engineering areas including cooperative control of unmanned aerial
vehicles (UAVs), scheduling of automated highway systems,
coordination/formation of underwater vehicles, attitude alignment
of satellite clusters and congestion control of communication
networks \cite{F. Giulietti}--\cite{R. Olfati-Saber 1}. There has
been considerable effort in modelling and exploring the collective
dynamics, and trying to understand how a group of autonomous
creatures or man-made mobile autonomous agents/robots can cluster
in formations without centralized coordination and control
\cite{C. W. Reynolds}--\cite{T. Chu3}.

In order to generate computer animation of the motion of flocks,
Reynolds \cite{C. W. Reynolds} modelled a flying bird as an object
moving in three dimensional environment based on the positions and
velocities of its nearby flockmates, and introduced the following
three rules (named steering forces) \cite{C. W. Reynolds}:

1) Collision Avoidance: avoid collisions with nearby flockmates,

2) Velocity Matching: attempt to match velocity with nearby
flockmates, and

3) Flock Centering: attempt to stay close to nearby flcokmates.\\
Subsequently, Vicsek {\it et al.} \cite{T. Vicsek} proposed a
simple model of autonomous agents (i.e., points or particles). In
the model, all agents move at a constant identical speed and each
agent updates its heading as the average of the heading of itself
with its nearest neighbors plus some additive noise. Moreover, the
authors used numerical simulations to demonstrate that all agents
eventually moved in the same direction, despite the absence of
centralized coordination and control. In fact, Vicsek's model can
be viewed as a special case of Reynolds's model, since it only
considers the regulation of velocity matching. Jadbabie {\it et
al.} \cite{A. Jadbabaie} and Savkin \cite{A. V. Savkin} used two
different methods to provide theoretical explanations for the
observed behaviors in Vicsek's model. Stimulated by the simulation
results in \cite{C. W. Reynolds}, Tanner {\it et al.} \cite{tanner
1}--\cite{tanner 2} considered a group of mobile agents moving on
the plane with double integrator dynamics. They introduced a set
of control laws that enable the group to generate stable flocking
motion, and provided theoretical justification. Nevertheless, it
is perhaps more reasonable and realistic to take the agents'
masses into account and consider the point mass model in which
each agent moves in high-dimensional space based on the Newton's
law. From \cite{H. Shi}, it is easy to see that, by using the
control laws given in \cite{tanner 1}, the group's final velocity
relies solely on the initial velocities of all agents in the
group. This means that these control laws cannot regulate the
final speed and heading of the group. On the other hand, in
reality, the motion of the group sometimes is inevitably
influenced by some external factors. Hence, it is not enough to
consider only the interactions among agents. In some cases, the
regulation of agents has certain purposes such as achieving
desired common speed and heading or arriving at a desired
destination. Therefore, the cooperation/coordination of multiple
mobile agents with some virtual leaders is an interesting and
important topic. There have been some papers dealing with this
issue in the literature. For example, Leonard and Fiorelli
\cite{N. E. Leonard and E. Fiorelli} viewed reference points as
virtual leaders used to manipulate the geometry of autonomous
vehicle group and direct the motion of the group. \cite{Y. Liu 2}
and \cite{A. Jadbabaie} considered the cohesion/coordination of a
group of mobile autonomous agents following an actual leader by
the so-called nearest neighbor rules.

In this paper, we investigate the collective behavior of
multi-agent systems in $n$-dimensional space with point mass
dynamics. By viewing the external control signals (or ``mission")
as virtual leaders, we show that all agents eventually move ahead
at a desired common velocity and maintain constant distances
between them. During the course of motion, each agent is
influenced by the external control signal and the motion of other
agents in the group. In order to generate the desired stable
flocking, we introduce a set of control laws such that each agent
regulates its position and orientation based on the desired
velocity and the information of a fixed set of ``neighbors". The
control laws are a combination of attractive/repulsive and
alignment forces. By using the control laws, all agent velocities
asymptotically approach the desired value, collisions can be
avoided between agents, and the final tight formation minimizes
all agent potentials. One salient feature of this paper is that
the self-organized global behavior is achieved via local feedback,
i.e. the desired emergent dynamics is produced through local
interactions and information exchange between the dynamic agents.

This paper is organized as follows: In Section 2, we formulate the
problem to be investigated. Some basic concepts and results in
graph theory are provided in Section 3. By using some specific
control laws, we analyze the system stability, the motion of the
center of mass (CoM), and the convergence rate of the system in
Section 4. We present some different control laws that can also
generate the desired stable flocking motion in Section 5. For the
case that not all agents know the desired velocity, we introduce a
set of control laws and study the system stability in Section 6.
Some numerical simulations are presented in Section 7. Finally, we
briefly summarize our results in Section 8.

\section{Problem Formulation}

We consider a group of $N$ agents moving in an $n$-dimensional
Euclidean space, each has point mass dynamics described by
\begin{equation}
\begin{array}{l}
\quad\,\dot{x}^i=v^i, \\
m_i\dot{v}^i=u^i,\ \ i=1, \cdots, N, \\
\end{array}
\label{eq1}
\end{equation}
where $x^i=(x_{1}^i, \cdots, x_{n}^i)^T\in R^n$ is the position
vector of agent $i$, $v^i=(v_{1}^i, \cdots, v_{n}^i)^T\in R^n$ is
its velocity vector, $m_i>0$ is its mass, and $u^i=(u_{1}^i,
\cdots, u_{n}^i)^T\in R^n$ is the (force) control input acting on
agent $i$. $x^{ij}=x^i-x^j$ denotes the relative position vector
between agents $i$ and $j$.

Our objective is to make the entire group move at a desired
velocity and maintain constant distances between the agents.
Additionally, we choose the control laws such that, during the
course of motion, collisions can be avoided between agents, and
the group final configuration minimizes all agent potentials. In
what follows, we will investigate the motion of the agent group in
two different cases, that is, we consider the group motion in
ideal case (i.e., velocity damping is ignored) and nonideal case,
respectively. For the two different cases, we propose two
different control laws such that the entire group moves at a
desired common velocity, and at the same time, collision-free can
be ensured between agents, and the group final configuration
minimizes all agent potentials.

We first consider the ideal case, that is, we ignore the velocity
damping. In this case, in order to achieve our objective, we try
to regulate each agent velocity to the desired velocity, reduce
the velocity differences between agents, and at the same time,
regulate their distances such that their potentials become
minimum. Hence, we choose the control law $u^i$ for agent $i$ to
be
\begin{equation}
u^i=\alpha^i+\beta^i+\gamma^i, \label{eq2}
\end{equation}
where $\alpha^i$ is used to regulate the potentials among agents,
$\beta^i$ is used to regulate the velocity of agent $i$ to the
weighted average of its flockmates, and $\gamma^i$ is used to
regulate the momentum of agent $i$ to the desired final momentum
(all to be designed later). $\alpha^i$ is derived from the social
potential fields which is described by artificial social potential
function $V^{i}$, a function of the relative distances between
agent $i$ and its flockmates. Collision-free and cohesion in the
group can be guaranteed by this term. $\beta^i$ reflects the
alignment or velocity matching with neighbors among agents.
$\gamma^i$ is designed to regulate the momentum among agents based
on the external signal (the desired velocity). By using such a of
momentum regulation, we can obtain the explicit convergence rate
of the CoM of the system.

\textbf{Remark 1}: The design of $\alpha^i$ and $\beta^i$
indicates that, during the course of motion, agent $i$ is
influenced only by its ``neighbors", whereas $\gamma^i$ reflects
the influence of the external signal on the agent motion.

Certainly, in some cases, the velocity damping can not be ignored.
For example, objects moving in viscous environment and mobile
objects with high speeds such as supersonic aerial vehicles, are
subjected to the influence of velocity damping. Then, in this
case, the model should be in the following form
\begin{equation}
\begin{array}{lcl}
\quad\,\dot{x}^i=v^i, \\
m_i\dot{v}^i=u^i-k_iv^i,\\
\end{array}
\label{eq3}
\end{equation}
where $k_i>0$ is the ``velocity damping gain", $-k_iv^i$ is the
velocity damping term, and $u^i$ is the control input for agent
$i$. Here we assume that the damping force is in proportion to the
magnitude of velocity. Moreover, since the ``velocity damping
gain" is determined by the shape and size of the object, the
property of the medium and some other factors, we assume that the
damping gains $k_i$, $i=1, \cdots, N$ are not equal to each other.
In order to achieve our objective, we need to compensate for the
velocity damping. Hence, we modify the control law $u^i$ to be
\begin{equation}
u^i=\alpha^i+\beta^i+\gamma^i+k_iv^i. \label{eq4}
\end{equation}

\section{Main Results}

In this section, we investigate the stability properties of
multiple mobile agents with point mass dynamics described in
(\ref{eq1}). We will present explicit control input in (\ref{eq2})
for the terms $\alpha^i$, $\beta^i$, and $\gamma^i$. We will
employ algebra and graph theory as basic tools for our discussion.
Some concepts and results in graph theory are given in the
Appendix.

Following \cite{tanner 1}, we make the following definitions and
assumptions.

\textbf{Definition 1} \cite{tanner 1}: (Neighboring graph) The
neighboring graph, $\cal{G}=(\cal{V}, \cal{E})$, is an undirected
graph consisting of

$i$) a set of vertices, ${\cal{V}}=\{n_1, \cdots, n_N\}$, indexed
by the agents in the group, and

$ii$) a set of edges, ${\cal{E}}=\{(n_i, n_j)\in{\cal{V}\times
\cal{V}} | n_j\sim n_i\}$, containing unordered pairs of vertices
that represent the neighboring relations.

In this paper, we consider a group of mobile agents with fixed
topology. We assume that the neighboring graph $\cal{G}$ is
connected, and hence does not change with time. Denote the set
${\cal{N}}_i\triangleq \{j | j\sim i \}\subseteq \{1, \cdots,
N\}\backslash \{i\}$ which contains all neighbors of agent $i$. If
agent $j$ is not a neighbor of agent $i$, we denote $j\nsim i$.

\textbf{Definition 2} \cite{tanner 1}: (Potential function)
Potential $V^{ij}$ is a differentiable, nonnegative, radially
unbounded function of the distance $\|x^{ij}\|$ between agents $i$
and $j$, such that

$i$) $V^{ij}(\|x^{ij}\|)\rightarrow \infty$ as
$\|x^{ij}\|\rightarrow 0$,

$ii$) $V^{ij}$ attains its unique minimum when agents $i$ and $j$
are located at a desired distance.

Functions $V^{ij}$, $i, j=1, \cdots, N$ are the artificial social
potential functions that govern the interindividual interactions.
Cohesion and separation can be achieved by artificial potential
fields \cite{E. Rimon}. In fact, cohesion can be ensured by the
connectivity of the neighboring graph, but collision-free can only
be guaranteed between interconnected agents. Collision can be
avoided between all agents only when the neighboring graph is
complete.

By the definition of $V^{ij}$, the total potential of agent $i$
can be expressed as
\begin{equation}
V^i=\sum_{j\in {\cal{N}}_i}V^{ij}(\|x^{ij}\|).
 \label{eq5}
\end{equation}

Agent dynamics in ideal case is different from that in nonideal
case, i.e., agents have different motion equations in the two
cases. Hence, in what follows, we will discuss the motion of the
agent group in two different cases, respectively.

Note that, in this section, we always assume that all agents can
receive the external signal, that is, they all know the desired
final velocity. In the case that not all agents know the mission,
we will discuss the flocking control problem in a separate
section.

\subsection{Ideal Case}

In this case, we take the control law $u^i$ for agent $i$ to be
\begin{equation}
u^i=-\sum_{j\in {\cal{N}}_i}w_{ij}(v^i-v^j)-\sum_{j\in
{\cal{N}}_i}\nabla_{x^i}V^{ij}-m_i(v^i-v^0),
 \label{eq6}
\end{equation}
where $v^0\in{R^n}$ is the desired common velocity and is a
constant vector, $w_{ij}\geq 0$, $w_{ij}=w_{ji}$, and $w_{ii}=0$,
$i, j=1, \cdots, N$ represent the interaction coefficients. And
$w_{ij}>0$ if agent $j$ is a neighbor of agent $i$, and is 0
otherwise. We denote $W=[w_{ij}]$. Thus, $W$ is symmetric, and by
the connectivity of the neighboring graph, $W$ is irreducible.

\subsubsection{Stability Analysis}

Before presenting the main results of this paper, we first prove
an important lemma.

\textbf{Lemma 1}: Let $A\in R^{n\times n}$ be any diagonal matrix
with positive diagonal entries. Then \[\left(A
{\mathrm{span}}\{\bf 1\}^\bot\right)\cap{\mathrm{span}}\{{\bf
1}\}=0,\] where ${\bf 1}=(1, \cdots, 1)^T\in{R^n}$,
${\mathrm{span}}\{{\bf 1}\}$ is the space spanned by the vector
${\bf 1}$, and ${\mathrm{span}}\{\bf 1\}^\bot$ is the orthogonal
complement space of ${\mathrm{span}}\{{\bf 1}\}$.

\textbf{Proof}: Let $p\in{\left(A {\mathrm{span}}\{\bf
1\}^\bot\right)\cap{\mathrm{span}}\{\bf 1\}}.$ Then
$p\in{\mathrm{span}\{\bf 1\}}$ and there is some $q\in{
{\mathrm{span}}\{\bf 1\}^\bot}$ such that $p=Aq$. It follows that
$q^TAq=q^Tp=0.$ Since $A$ is positive definite by assumption, we
have $q=0$ and hence $p=0$. \hfill $\square $

\textbf{Theorem 1}: By taking the control law in (\ref{eq6}), all
agent velocities in the group described in (\ref{eq1})
asymptotically approach the desired common velocity,
collision-free is ensured between neighboring agents, and the
group final configuration minimizes all agent potentials.

This theorem becomes apparently true after Theorem 2 is proved, so
we proceed to present Theorem 2.

We define the following error vectors:
\begin{equation}
e_{p}^i=x^i-v^0t, \label{eq7}
\end{equation}
\begin{equation}
e_{v}^i=v^i-v^0, \label{eq8}
\end{equation}
where $t$ is time variable and $v^0$ is the desired common
velocity. $e_{p}^i$ represents the relative position vector
between the actual position of agent $i$ and its desired position.
$e_{v}^i$ represents the velocity difference vector between the
actual velocity and the desired velocity of agent $i$. It is easy
to see that $\dot{e}_{p}^i=e_{v}^i$, and
$\dot{e}_{v}^i=\dot{v}^i$. Hence, the error dynamics is given by
\begin{equation}
\begin{array}{l}
\dot{e}_{p}^i=e_{v}^i, \\[6mm]
\dot{e}_{v}^i=\frac{1}{m_i}u^i, \ \ i=1, \cdots, N.\\
\end{array}
\label{eq9}
\end{equation}
By the definition of $V^{ij}$, it follows that
\[V^{ij}(\|x^{ij}\|)=V^{ij}(\|e_{p}^{ij}\|):=\widetilde{V}^{ij},\]
where $e_{p}^{ij}\triangleq e_{p}^i-e_{p}^j$, and hence
$\widetilde{V}^i=V^i$ and
$\nabla_{e_{p}^i}\widetilde{V}^{ij}=\nabla_{x^i}V^{ij}.$ Thus, the
control input for agent $i$ in the error system has the following
form
\begin{equation}
u^i=-\sum_{j\in {\cal{N}}_i}w_{ij}(e_{v}^i-e_{v}^j)-\sum_{j\in
{\cal{N}}_i}\nabla_{e_{p}^i}\widetilde{V}^{ij}-m_ie_{v}^i.
 \label{eq10}
\end{equation}

Consider the following positive semi-definite function
\begin{equation}
J=\frac{1}{2}\sum_{i=1}^N\left(\widetilde{V}^i+m_ie_{v}^{iT}e_{v}^i\right).
\label{eq11}
\end{equation}
It is easy to see that $J$ is the sum of the total artificial
potential energy and the total kinetic energy of all agents in the
error system. Define the level set of $J$ in the space of agent
velocities and relative distances in the error system
\begin{equation}
\Omega=\left\{(e_{v}^i, e_{p}^{ij}) | J\leq c \right\}.
 \label{eq12}
\end{equation}
In what follows, we will prove that the set $\Omega$ is compact.
In fact, the set $\{e_{v}^i, e_{p}^{ij}\}$ with $J\leq c$ ($c>0$)
is closed by continuity. Moreover, boundedness can be proved by
connectivity. More specifically, from $J\leq c$, we have
$\widetilde{V}^{ij}\leq c$. Moreover, since the potential function
$V^{ij}$ is radially unbounded, $\widetilde{V}^{ij}$ is also
radially unbounded, and there is a positive constant $d_{ij}$ such
that $\|e_{p}^{ij}\|\leq d_{ij}$. Denote
$d=\max\limits_{j\in{{\cal{N}}_i}}d_{ij}$. Since the neighboring
graph is connected, there must be a path connecting any two agents
$i$ and $j$, and its length does not exceed $N-1$. Hence, we have
$\|e_{p}^{ij}\|\leq (N-1)d$. By similar analysis, we have
$e_{v}^{iT}e_{v}^i\leq 2c/m_i$, thus
$\|e_{v}^i\|\leq{\sqrt{2c/m_i}}$.

By the symmetry of $\widetilde{V}^{ij}$ with respect to
$e_{p}^{ij}$ and by $e_{p}^{ij}=-e_{p}^{ji}$, it follows that
\begin{equation}
\frac{\partial \widetilde{V}^{ij}}{\partial
e_{p}^{ij}}=\frac{\partial \widetilde{V}^{ij}}{\partial
e_{p}^i}=-\frac{\partial \widetilde{V}^{ij}}{\partial e_{p}^j},
 \label{eq13}
\end{equation}
and therefore
\[\frac{d}{dt}\sum_{i=1}^N\frac{1}{2}\widetilde{V}^i=\sum_{i=1}^N\nabla_{e_{p}^i}\widetilde{V}^i\cdot
e_{v}^i.\]

\textbf{Theorem 2}: By taking the control law in (\ref{eq10}), all
agent velocities in the system described in (\ref{eq9})
asymptotically approach zero, collision-free is ensured between
neighboring agents, and the group final configuration minimizes
all agent potentials.

\textbf{Proof}: Choosing the positive semi-definite function $J$
defined as in (\ref{eq11}) and calculating the time derivative of
$J$ along the solution of the error system (\ref{eq9}), we have
\begin{equation}
\begin{array}{rl}
\dot{J}=& \! \! \!
-\displaystyle\sum\limits_{i=1}^N\sum\limits_{j\in
{\cal{N}}_i}w_{ij}e_{v}^{iT}(e_{v}^i-e_{v}^j)-\sum_{i=1}^Nm_ie_{v}^{iT}e_{v}^i\\
=& \! \! \! \displaystyle-e_{v}^T(L\otimes
I_n)e_v-e_{v}^T(M\otimes I_n)e_v,
\end{array}
\label{eq14}
\end{equation}
where $e_v=(e_{v}^{1T}, \cdots, e_{v}^{NT})^T$ is the stack vector
of all agent velocity vectors in the error system; $L=[l_{ij}]$
with
\begin{equation}
\begin{array}{l}
{l_{ij}}=\left\{
\begin{array}{l}
-w_{ij}, \\
\sum_{k=1,k\neq i}^Nw_{ik},%
\end{array}
\begin{array}{l}
\;i\neq j, \\
\;i=j;
\end{array}
\right.  \\
\end{array}
\label{eq15}
\end{equation}
$M={\mathrm {diag}}(m_1, \cdots, m_N)$; $I_n$ is the identity
matrix of order $n$ and $\otimes$ stands for the Kronecker
product.

By the definition of matrix $L$, it is easy to see that $L$ is
symmetric, each row sum is equal to 0, the diagonal entries are
positive, and all the other entries are nonpositive. By matrix
theory \cite{R. Horn and C. R. Johnson}, all eigenvalues of $L$
are nonnegative. Hence, matrix $L$ is positive semi-definite. By
the connectivity of the neighboring graph and the symmetry of
matrix $L$, it follows that $L$ is irreducible and the eigenvector
associated with the single zero eigenvalue is ${\bf 1}_N$. On the
other hand, it is known that the identity matrix $I_n$ has an
eigenvalue $\mu=1$ of $n$ multiplicity and $n$ linearly
independent eigenvectors
\begin{equation*}
p^1=\left[
\begin{array}{cc}
1 \\
0\\
\vdots\\
0%
\end{array}%
\right] ,\;\;p^2=\left[
\begin{array}{c}
0 \\
1\\
\vdots\\
0%
\end{array}%
\right] ,\;\;\cdots,\;\; p^n=\left[
\begin{array}{cc}
0\\
\vdots\\
0\\
1%
\end{array}%
\right] .
\end{equation*}%
By matrix theory \cite{R. Horn and C. R. Johnson}, the eigenvalues
of $L\otimes I_n$ are nonnegative, $\lambda=0$ is an eigenvalue of
multiplicity $n$ and the associated eigenvectors are
\[q^1=[p^{1T}, \cdots, p^{1T}]^T, \cdots, q^n=[p^{nT}, \cdots,
p^{nT}]^T.\] Furthermore, it is easy to see that matrix $M$ is
positive definite and hence $-e_{v}^T(M\otimes I_n)e_v\leq 0$.
Thus $\dot{J}\leq 0$, and $\dot{J}=0$ implies that
$e_{v}^1=e_{v}^2=\cdots=e_{v}^N$ and they all must equal zero.
This occurs only when $v^1=v^2=\cdots=v^N=v^0$, that is, the
vector $e_{vk}=(e_{vk}^1, \cdots, e_{vk}^N)$ $(k=1, \cdots, n)$,
which is composed of all the corresponding $k$th components
$e_{vk}^1, \cdots, e_{vk}^N$ of $e_{v}^1, \cdots, e_{v}^N$, is
contained in span\{{\bf 1}\}, where ${\bf 1}=(1, \cdots, 1)^T\in
R^N$ and each entry $e_{vk}^i$ of $e_{vk}$ equals zero. It follows
that $\dot{e}_{p}^{ij}=0$, $\forall (i, j)\in{N\times N}$.

We use LaSalle's invariance principle \cite{H. K. Khalil} to
establish convergence of the system trajectories to the largest
positively invariant subset of the set defined by $E=\{e_v |
\dot{J}=0\}$. In $E$, the agent velocity dynamics in the error
system is
\[\dot{e}_{v}^i=\frac{1}{m_i}u^i=-\frac{1}{m_i}\sum_{j\in
{\cal{N}}_i}\nabla_{e_{p}^i}\widetilde{V}^{ij}=-\frac{1}{m_i}\nabla_{e_{p}^i}\widetilde{V}^i\]
and therefore it follows that
\begin{equation}
\begin{array}{rl}
\dot{e}_v=-((M^{-1}B)\otimes I_n)\left[
\begin{array}{c}
\vdots \\
\nabla_{e_{p}^{ij}}\widetilde{V}^{ij}\\
\vdots
\end{array}%
\right],
\end{array} \label{eq16}
\end{equation}
where $M^{-1}={\mathrm{diag}} (\frac{1}{m_1}, \cdots,
\frac{1}{m_N})$ is the inverse of matrix $M$, and matrix $B$ is
the incidence matrix of the neighboring graph. Hence
\[\dot{e}_{vk}=-(M^{-1}B)[\nabla_{e^{ij}}\widetilde{V}^{ij}]_k,\ \  k=1, \cdots, n.\]
Thus, $\dot{e}_{vk}\in {\mathrm{range}}(M^{-1}B)$, $k=1, \cdots,
n$. By matrix theory and by the connectivity of the neighboring
graph $\cal{G}$, we have
\[{\mathrm{range}}(M^{-1}B)=M^{-1}{\mathrm{range}}B=M^{-1}{\mathrm{range}}(BB^T)=M^{-1}{\mathrm{span}}\{{\bf
1}\}^\bot\] and therefore
\begin{equation}
\dot{e}_{vk}\in M^{-1}{\mathrm{span}}\{{\bf 1}\}^\bot, \ \ k=1,
\cdots, n. \label{eq17}
\end{equation}
In any invariant set of $E$, by $e_{vk}\in {\mathrm{span}}\{{\bf
1}\}$, we have
\begin{equation}
\dot{e}_{vk}\in {\mathrm{span}}\{{\bf 1}\}. \label{eq18}
\end{equation}
Furthermore, by Lemma 1, we get from (\ref{eq17}) and (\ref{eq18})
that
\[\dot{e}_{vk}\in (M^{-1}{\mathrm{span}}\{{\bf 1}\}^\bot) \cap
 {\mathrm{span}}\{{\bf 1}\}\equiv {\bf 0}, \ \ k=1, \cdots, n.\]
Thus, in steady state, all agent velocities in the error system no
longer change and equal zero, and moreover, from (\ref{eq16}), the
potential $\widetilde{V}^i$ of each agent $i$ is minimized.
Collision-free can be ensured between neighboring agents since
otherwise it will result in $\widetilde{V}^i\rightarrow\infty$.
\hfill $\square $

From the proof of Theorem 2, it follows that, in steady state, all
agent actual velocities no longer change and are equal to the
desired velocity.

\textbf{Remark 2}: Only when the neighboring graph is complete,
collision avoidance between all agents can be guaranteed with the
control laws above.

\textbf{Remark 3}: If we take the control law for agent $i$ to be

\begin{equation}
u^i=-\sum_{j\in {\cal{N}}_i}(v^i-v^j)-\sum_{j\in
{\cal{N}}_i}\nabla_{x^i}V^{ij}-m_i(v^i-v^0),\label{eq19}
\end{equation}
we can also get the same conclusion as in Theorem 1. Here, we
still consider the error system (\ref{eq9}). In fact, if we take
the same Laypunov function $J$ defined as in Theorem 2 and take
the control law in (\ref{eq19}), we obtain that
\[\dot{J}=-e_{v}^T(\overline{L}\otimes I_n)e_v-e_{v}^T(M\otimes
I_n)e_v,\] where $\overline{L}=\Delta-{\cal{A}}$ is the Laplacian
matrix of the neighboring graph, $\Delta={\mathrm{diag}}
(N_1,\cdots,N_N)$, $N_i$ denotes the valance of vertex $i$ in the
graph, and $\cal{A}$ is the adjacency matrix of the graph. Using a
similar analysis method as in Theorem 2, we can obtain the same
conclusion of stable flocking. Note that, in comparison with
Theorem 2, we have the following difference on the decaying rates
of the energy function $J$
\[-e_{v}^T(L\otimes I_n)e_v+e_{v}^T(\overline{L}\otimes
I_n)e_v=-e_{v}^T(\widetilde{L}\otimes I_n)e_v,\] where
$\widetilde{L}=[\widetilde{l}_{ij}]$ with
\begin{equation*}
\begin{array}{lcl}
&  & {\widetilde{l}_{ij}}=\left\{
\begin{array}{l}
-w_{ij}+1, \\
0\\
\sum_{k=1,k\neq i}^Nw_{ik}-N_i,%
\end{array}
\begin{array}{l}
\;i\neq j\ {\mathrm{and}}\ j\sim i, \\
\;i\neq j\ {\mathrm{and}}\ j\nsim i, \\
\;i=j.
\end{array}
\right.  \\
\end{array}
\end{equation*}
It is easy to see that, by using the different control laws in
(\ref{eq6}) and (\ref{eq19}), the decaying rates of the total
energy $J$ may be different. Hence, the interaction coefficients
$w_{ij}$ can influence the convergence rate of the system.

\subsubsection{The Motion of the Center of Mass}

In what follows, we will analyze the motion of the center of mass
(CoM) of system (\ref{eq1}).

The position vector of the CoM in system (\ref{eq1}) is defined as
\[x^*=\frac{\sum_{i=1}^Nm_ix^i}{\sum_{i=1}^Nm_i}.\] Thus, the
velocity vector of the CoM is
\[v^*=\frac{\sum_{i=1}^Nm_iv^i}{\sum_{i=1}^Nm_i}.\] By using
control law (\ref{eq6}), we obtain
\begin{equation*}
\begin{array}{rl}
\dot{v}^*=
\displaystyle\frac{-1}{\big(\sum_{i=1}^Nm_i\big)}\sum_{i=1}^N\left[\sum_{j\in
{\cal{N}}_i}w_{ij}(v^i-v^j)+\sum_{j\in
{\cal{N}}_i}\nabla_{x^i}V^{ij}+m_i(v^i-v^0)\right].%
\end{array}
\end{equation*}%
By the symmetry of matrix $W$ and the symmetry of function
$V^{ij}$ with respect to $x^{ij}$, we get
\begin{equation}
\dot{v}^*=-v^*+v^0. \label{eq20}
\end{equation}
Suppose the initial time $t_0=0$, and $v^*(0)=v_{0}^*$. By solving
(\ref{eq20}), we get
\[v^*=v^0+(v_{0}^*-v^0)e^{-t}.\]
Thus, it follows that, if $v_{0}^*=v^0$, then the velocity of the
CoM is invariant and equals $v^0$ for all the time; if
$v_{0}^*\neq v^0$, then the velocity of the CoM exponentially
converges to the desired velocity $v^0$ with convergence exponent
1.

Therefore, from the analysis above, we have the following theorem.

\textbf{Theorem 3}: By taking the control law in (\ref{eq6}), if
the initial velocity of the CoM is equal to the desired velocity,
then it is invariant for all the time; otherwise it will
exponentially converge to the desired velocity.

\textbf{Remark 4}: Note that, by the calculation above, we can see
that the velocity variation of the CoM does not rely on the
neighboring relations or the magnitudes of the interaction
coefficients. Even if the neighboring graph is not connected, the
velocity of the CoM still equals the desired velocity or
exponentially converges to it, and the final velocities of all
connected agent groups equal the desired velocity as well.
However, in this case, the distance between disconnected subgroups
might be very far.

\subsubsection{Convergence Rate Analysis}

From the discussion above, we know that the coupling coefficients
can influence the convergence rate of system (\ref{eq1}). In what
follows, we will present some qualitative analysis of the
influence of the weights on the convergence rate of the system.

Let us again consider the dynamics of the error system. From the
analysis in Theorem 2, we know that $\dot{J}\leq 0$, and
$\dot{J}=0$ occurs only when $e_{v}^1=e_{v}^2=\cdots=e_{v}^N=0$,
that is, only when all agents have reached the desired velocity.
In other words, if there exists one agent whose velocity is
different from the desired velocity, then the energy function $J$
is strictly monotone decreasing with time. Of course, before the
group forms the final tight configuration, there might be the case
that all agent velocities have reached the desired value, but due
to the regulation of the potentials among neighboring agents, it
instantly changes into the case that not all agents have the
desired velocity. Hence, the decaying rate of energy is equivalent
to the convergence rate of the system. It is easy to see that,
when not all agents have reached the desired velocity, for any
solution of the error system, $e_v$ must be in the subspace
spanned by the eigenvectors of $L\otimes I_n$ corresponding to the
nonzero eigenvalues. Thus, from (\ref{eq14}), we have $\dot{J}\leq
-\lambda_2e_{v}^Te_v-e_{v}^T(M\otimes I_n)e_v $, where $\lambda_2$
denotes the second smallest real eigenvalue of matrix $L$.
Therefore, we have the following conclusion: The convergence rate
of the system relies on the second smallest real eigenvalue of
matrix $L$ defined as in (\ref{eq15}) as well as agent masses, and
it is always not faster than the convergence rate of the CoM.
Furthermore, if the initial velocity of the CoM is not equal to
the desired velocity, then the fastest convergence rate of the
system does not exceed the exponential convergence rate with
convergence exponent 1.

\textbf{Remark 5}: Note that when the group has achieved the final
steady state, the control input above equals zero.

\subsection{ Nonideal Case}

Sometimes, the velocity damping should not be ignored. Then, in
this case, in order to make the group generate the desired stable
flocking motion, the velocity damping need to be cancelled by some
terms in the control laws. Hence, we modify the control law as in
(\ref{eq4}), where $\alpha^i$, $\beta^i$, and $\gamma^i$ are
defined as in (\ref{eq6}), that is, the control law acting on
agent $i$ is
\begin{equation}
u^i=-\sum_{j\in {\cal{N}}_i}w_{ij}(v^i-v^j)-\sum_{j\in
{\cal{N}}_i}\nabla_{x^i}V^{ij}-m_i(v^i-v^0)+k_iv^i. \label{eq21}
\end{equation}
Then, the total force acting on agent $i$ is
\[
u^i=-\sum_{j\in {\cal{N}}_i}w_{ij}(v^i-v^j)-\sum_{j\in
{\cal{N}}_i}\nabla_{x^i}V^{ij}-m_i(v^i-v^0).
\]

All the results in ideal case can be analogously extended to the
nonideal case. Namely, following Theorems 1 and 2, we can easily
obtain the desired stable flocking motion, that is, when the
velocity damping is taken into account, by using control law
(\ref{eq21}), all agent velocities in the group described in
(\ref{eq3}) asymptotically approach the desired value,
collision-free can be ensured between neighboring agents, and the
group final configuration minimizes all agent potentials.
Furthermore, following Theorem 3 and the convergence rate analysis
above, we conclude that the convergence rate of the system relies
on the interaction coefficients and agent masses, and when the
initial velocity of the CoM is not equal to the desired velocity,
the fastest convergence rate of the system does not exceed the
exponential convergence rate with convergence exponent 1.

Note also that, because the velocity damping is concelled by some
terms in the control law, the velocity damping cannot influence
the convergence rate of system (\ref{eq3}).

\textbf{Remark 5}: In steady state, the group keeps on moving at a
desired velocity. During this period, the control laws' role is
only to cancel the velocity damping.

\section{Discussions on Various Control Laws}

In the sections above, we introduced a set of control laws that
enable the group to generate the desired stable flocking motion.
However, it should be clear that control law (\ref{eq6}) is not
the unique control law to produce the desired motion for the
group. In this section, we provide some more useful control laws.
For simplicity, we only present the control laws for the group
moving in the ideal case, since in the nonideal case, we only need
to add the terms $k_iv^i$ $(i=1, \cdots, N)$ to cancel the
velocity damping.

In the sequel, we will propose three different control laws that
can achieve our control objective. The analysis and proofs are
quite similar for these control laws, so we only present the
control laws and their corresponding Lyapunov functions.

1) In the control laws above, $\gamma^i$ is used to regulate the
momentum of agent $i$. However, we can also use $\gamma^i$ to
directly regulate the velocity of agent $i$ to the desired value.
Hence, we take the control law acting on agent $i$ to be
\begin{equation}
u^i=-\sum_{j\in {\cal{N}}_i}w_{ij}(v^i-v^j)-\sum_{j\in
{\cal{N}}_i}\nabla_{x^i}V^{ij}-(v^i-v^0). \label{eq22}
\end{equation}

We still consider the error system (\ref{eq9}) and choose Lyapunov
function (\ref{eq11}). By similar calculation, we get
\[\dot{J}=-e_{v}^T(L\otimes I_n)e_v-e_{v}^Te_v.\] Using the same
analysis method as in Theorem 2, we obtain that $\dot{J}\leq 0$,
and $\dot{J}=0$ implies that $e_{v}^1=e_{v}^2=\cdots=e_{v}^N=0$.
The rest analysis is similar to Theorem 2, thus is omitted.

\textbf{Remark 6}: Note that, the control law in (\ref{eq22}) can
make the group generate the desired stable flocking motion. But we
cannot explicitly estimate the convergence rate of the CoM by
using this control law.

2) Suppose that $\alpha^i$ and $\beta^i$ rely on agent $i$'s mass.
The control law acting on agent $i$ has the following form
\begin{equation}
 u^i=-\sum_{j\in {\cal{N}}_i}m_iw_{ij}(v^i-v^j)-\sum_{j\in
{\cal{N}}_i}m_i\nabla_{x^i}V^{ij}-m_i(v^i-v^0). \label{eq23}
\end{equation}

In this case, for the error system (\ref{eq9}), we choose the
following Lyapunov function
\begin{equation}
J=\frac{1}{2}\sum_{i=1}^N\left(\widetilde{V}^i+e_{v}^{iT}e_{v}^i\right).
\label{eq24}
\end{equation}
By similar calculation, we have \[\dot{J}=-e_{v}^T(L\otimes
I_n)e_v-e_{v}^Te_v.\] Following the analysis method in Theorem 2,
we can show that the desired stable flocking motion will be
achieved.

\textbf{Definition 3}: Define the center of the system of agents
as $\overline{x}=(\sum_{i=1}^Nx^i)/N.$

\textbf{Definition 4}: The average velocity of all agents is
defined as $\overline{v}=(\sum_{i=1}^Nv^i)/N.$

It is obvious that the velocity of the system center is just the
average velocity of all agents.

Using the control law in (\ref{eq23}), we have
$\dot{\overline{v}}=-\overline{v}+v^0$. Suppose the initial time
$t_0=0$ and $\overline{v}(0)=\overline{v}_0$. We get
\[\overline{v}=v^0+(\overline{v}_0-v^0)e^{-t}.\]

It is obvious that, if $\overline{v}_0=v^0$, then the velocity of
the system center is equal to the desired velocity $v^0$ for all
the time, and if $\overline{v}_0\neq v^0$, then the velocity of
the system center exponentially converges to the desired velocity
with convergence exponent 1.

3) Suppose that $\alpha^i$ and $\beta^i$ rely on agent $i$'s mass,
and $\gamma^i$ is used to regulate the velocity of agent $i$ to
the desired velocity. The control law $u^i$ is then taken to be
\begin{equation}
u^i=-\sum_{j\in {\cal{N}}_i}m_iw_{ij}(v^i-v^j)-\sum_{j\in
{\cal{N}}_i}m_i\nabla_{x^i}V^{ij}-(v^i-v^0).\label{eq25}
\end{equation}

We consider the error system (\ref{eq9}) and choose the
corresponding Lyapunov function (\ref{eq24}). Then,
\[\dot{J}=-e_{v}^T(L\otimes I_n)e_v-e_{v}^T(M^{-1}\otimes
I_n)e_v,\] where $M^{-1}$ is the inverse of matrix $M$. The rest
analysis is similar, thus is omitted.

\textbf{Remark 7}: Note that, using the control law in
(\ref{eq25}), the convergence rates of the CoM and the system
center both cannot be explicitly estimated.

From the analysis above, we conclude that the control law in
(\ref{eq6}) is the best one among the various control laws. On the
one hand, control law (\ref{eq6}) can be given certain physical
explanations, on the other hand, the corresponding Lyapunov
function has certain physical meaning. More importantly, by using
the control law in (\ref{eq6}), the convergence rate of the CoM of
the system can be accurately estimated.

\section{Extensions and Discussions}

In this section, we investigate the case that not all agents know
the desired velocity. We assume that the neighboring graph is
connected and the neighboring relations are fixed.

We first divide the group into two subgroups. Subgroup One
consists of all agents that can detect the reference signal, i.e.,
all agents who know the desired velocity belong to Subgroup One.
Subgroup Two contains all agents that can not detect the reference
signal. Hence, each agent in Subgroup One regulates its state
based on the reference signal and the information of its
``neighbors", whereas each agent in Subgroup Two regulates its
state only based on its ``neighbors".

We assume that there exists at least one agent who knows the
desired velocity. In the case that there is no external signal
acting on the group, the collective dynamic behaviors of the agent
group have been analyzed in \cite{H. Shi}.

Without loss of generality, suppose that agent $i$ $(i=1, \cdots,
N_1)$ $(1 \leq N_1 < N)$ are contained in Subgroup One, and agent
$j$ $(j=N_1+1, \cdots, N)$ are contained in Subgroup Two. The
control law acting on agent $i$ in Subgroup One is taken to be
\[ u^i=-\sum_{k\in
{\cal{N}}_i}w_{ik}(v^i-v^k)-\sum_{k\in
{\cal{N}}_i}\nabla_{x^i}V^{ik}-m_i(v^i-v^0),\] and the control law
acting on agent $j$ in Subgroup Two is taken to be
\[ u^j=-\sum_{k\in
{\cal{N}}_j}w_{jk}(v^j-v^k)-\sum_{k\in
{\cal{N}}_j}\nabla_{x^j}V^{jk}.\] We can lump the two equations
above into one
\begin{equation}
u^i=-\sum_{j\in {\cal{N}}_i}w_{ij}(v^i-v^j)-\sum_{j\in
{\cal{N}}_i}\nabla_{x^i}V^{ij}-h_im_i(v^i-v^0) \label{eq26}
\end{equation}
for $i=1, \cdots, N$, where $h_i$ is defined as
\begin{equation*}
\begin{array}{l}
h_i=\left\{
\begin{array}{l}
1, \\
0,%
\end{array}
\begin{array}{l}
\;{\mathrm{if\, agent}}\, i \,{\mathrm{\, is\, contained\, in\, Subgroup\, One}}, \\
\;{\mathrm{if\, agent}}\, i \,{\mathrm{\, is\, contained\, in\,
Subgroup\, Two}}.
\end{array}
\right.  \\
\end{array}
\end{equation*}%

We still consider the error system (\ref{eq9}). Using control law
(\ref{eq26}) and taking Lyapunov function (\ref{eq11}), we have
\[\dot{J}=-e_{v}^T(L\otimes I_n)e_v-e_{v}^T(\widehat{M}\otimes
I_n)e_v,\] where $L, e_{v}^i$, and $I_n$ are defined as before,
$\widehat{M}={\mathrm{diag}}(h_1m_1, \cdots, h_Nm_N)$. By the
definition of $h_i$ and $m_i$, it follows that matrix
$\widehat{M}$ is positive semi-definite. Following similar
analysis as in the previous sections, we can conclude that the
desired stable flocking motion can be achieved.

\textbf{Remark 8}: If there exists only one agent in the group who
can detect the external reference signal, the group can still
generate the desired stable flocking motion. This is of practical
interest in control of multi-agent systems.

\textbf{Remark 9}: Even if only one agent in the group cannot
detect the external reference signal, it is difficult to
explicitly estimate the convergence rate of the CoM.

It should be noted that there is no actual leader among agents,
all agents play the same role. However, we can view the external
reference signal as a virtual leader.

The results in this section suggest that, if we want to control a
group of mobile agents to move at a given velocity, we only need
to send our mission signal to any one of them. Then the signal can
be propagated through the neighboring interactions. This is of
practical interest in control of multiple mobile robots or a large
population of animals (think how do you pass through a crowds of
people? and how the shepherding dog steer a large group of sheep
back home?).

\section{Simulations}

In this section, we will present some numerical simulations for
the system described by (\ref{eq1}) in order to illustrate the
theoretic results obtained in the previous sections.

These simulations are performed with ten agents moving on the
plane whose initial positions, velocities and the neighboring
relations are set randomly, but they satisfy: 1) all initial
positions are set within a ball of radius $R=15$[m] centered at
the origin, 2) all initial velocities are set with arbitrary
directions and magnitudes within the range of (0, 10)[m/s], and 3)
the neighboring graph is connected. All agents have different
masses and they are set randomly in the range of (0, 1)[kg].

The following simulations are all performed with the same group,
and the group has the same initial state, including all agent
initial positions, velocities, and the fixed neighboring relations
between agents. However, different control laws are taken in the
form of (\ref{eq6}) or (\ref{eq26}) with the explicit potential
function
\[V^{ij}=5\ln{\|x^{ij}\|^2}+\frac{5}{\|x^{ij}\|^2},\ \ i, j=1,
\cdots, 10.\]

The interaction coefficient matrix $W$ is generated randomly such
that $w_{ii}=0$, $w_{ij}=w_{ji}$, and the nonzero $w_{ij}$ satisfy
$0<w_{ij}<1$ for all $i, j=1, \cdots, 10$. We run all simulations
for 250 seconds.

Fig. 1 presents the group initial state including the initial
positions and velocities of all agents, and the neighboring
relations between agents. Figs. 2--4, 5--7 and 8--10 show the
motion trajectories of all agents, the final configurations of the
group, and all agent velocities in three different simulations,
respectively. Figs. 11 and 12 depict the motion trajectories of
the CoM in the three simulations, whereas Figs. 13 and 14 are the
velocity curves of the CoM. Note that, in the velocity curve
figures, the solid arrow indicates the tendency of velocity
variation, and the meanings of the other arrows, dashed lines, and
solid lines are all presented in the figures.

In Fig. 1, the solid lines represent the neighboring relations and
the dotted arrows represent the initial velocity vectors. Figs.
2--4 describe the group state in the case that the motion of the
group is not influenced by any external signal and only relies on
the interactions between agents. It can be seen from them that,
during the course of motion, all agents regulate their positions
to minimize their potentials, regulate their velocities to reduce
the differences , and move ahead with a steady state
configuration. Moreover, the final common velocity is equal to the
initial velocity of the CoM of the system.

When we send a signal to the group and try to make all agents move
at a desired velocity, Figs. 5--7 show the results in our
simulation with the control laws taken in the form of (\ref{eq6}),
whereas Figs. 8--10 show the simulation results with the control
laws taken in the form of (\ref{eq26}) and with the assumption
that there is only one agent who knows the desired velocity. It
can be seen from them that all agents regulate their positions to
minimize their potentials and eventually move ahead with a steady
state configuration. Figs. 7 and 10 are the velocity curves, and
they distinctly demonstrate that all agent velocities
asymptotically approach the desired velocity.

Note that the final configurations of the group are different in
the three simulations. This is because, during the course of
motion, each agent regulates its position only based on the
information of its ``neighbors" in the group, hence collisions
cannot be avoided between the agents having no neighboring
relations.

Fig. 11 shows the motion trajectories of the CoM in the
simulations where the star represents the initial position of the
CoM, and Fig. 12 is the magnification of the trajectories of the
CoM at the initial time. Fig. 13 shows the velocity curves of the
CoM where the star represents the initial velocity of the CoM, and
Fig. 14 is the magnification of the velocity curves of the CoM at
the initial time. In these four figures, (a), (b), and (c)
represent the corresponding states of the CoM in the three
simulations, respectively. It can be seen from them that, when
there is no external signal acting on the group, the velocity of
the CoM is always invariant and is equal to the final common
velocity, otherwise, the velocity of the CoM converges to the
desired velocity. Apparently, the convergence rate of the CoM is
faster than the convergence rate of the system.

Hence, numerical simulations also indicate that, by using the
control law in (\ref{eq6}), the desired stable flocking motion can
be achieved.

\section{Conclusions}

We have investigated the collective behavior of multiple dynamic
agents moving in high-dimensional space with point mass dynamics,
and presented some control laws which ensure the group to generate
the desired stable flocking motion. The group dynamic properties
are characterized in two different cases. When the velocity
damping is negligible, using a set of coordination control laws,
we can make the group generate the desired stable flocking motion.
The control laws are a combination of attractive/repulsive and
alignment forces, and they ensure that all agent velocities
asymptotically approach the desired velocity, collisions are
avoided between neighboring agents, and the final tight formation
minimizes all agent potentials. Moreover, we showed that, when the
initial velocity of the center of mass is not equal to the desired
velocity, it will exponentially converge to the desired velocity.
When the velocity damping is taken into account, we can properly
modify the control laws in order to generate the desired stable
flocking. Subsequently, we investigated the motion of the group in
the case that not all agents know the desired final velocity, and
showed that the desired stable flocking motion can still be
achieved by our control laws. Finally, numerical simulations were
worked out to further illustrate our theoretical results. Our
method is general, integrating both algebraic theory and graph
theory, and is applicable to dealing with more complex agent
dynamics, information topology and interaction mechanisms.

\section{Appendix: Graph Theory Preliminaries}

In this section, we briefly summarize some basic concepts and
results in graph theory that have been used in this paper. More
comprehensive discussions can be found in \cite{C. Godsil and G.
Royle}.

A {\it graph} $\cal{G}$ consists of a {\it vertex set}
${\cal{V}}=\{n_1, n_2, \cdots, n_m\}$ and an {\it edge set}
${\cal{E}}=\{(n_i, n_j): n_i, n_j \in {\cal{V}}\}$, where an {\it
edge} is an unordered pair of distinct vertices of $\cal{V}$. If
$n_i, n_j\in{\cal{V}}$, and $(n_i, n_j)\in{\cal{E}}$, then we say
that $n_i$ and $n_j$ are {\it adjacent} or that $n_j$ is a {\it
neighbor} of $n_i$, and denote this by writing $n_j\sim n_i$. A
graph is called {\it complete} if every pair of vertices are
adjacent. The {\it valence} of vertex $n_i$ of ${\cal{G}}$ is
defined as the number of edges of ${\cal{G}}$ which are incident
with $n_i$, where an edge is {\it incident} with vertex $n_i$ if
one of the two vertices of the edge is $n_i$. The {\it adjacency
matrix} of ${\cal{G}}$ is an $m\times m$ matrix whose $ij$th entry
is 1 if $(n_i, n_j)$ is one of ${\cal{G}}$'s edges and is 0 if it
is not. A {\it path of length} $r$ from $n_i$ to $n_j$ in a graph
is a sequence of $r+1$ distinct vertices starting with $n_i$ and
ending with $n_j$ such that consecutive vertices are adjacent. If
there exists a path between any two vertices of ${\cal {G}}$, then
${\cal {G}}$ is {\it connected}.

An {\it oriented graph} is a graph together with a particular
orientation, where the {\it orientation} of a graph ${\cal {G}}$
is the assignment of a direction to each edge, so edge $(n_i,
n_j)$ is an directed edge from $n_i$ to $n_j$. The {\it incidence
matrix} $B$ of an oriented graph $\cal{G}$ is the $\{0, \pm
1\}$-matrix with rows and columns indexed by the vertices and
edges of $\cal{G}$, respectively, such that the $ij$-entry is
equal to 1 if edge $j$ is ending on vertex $n_i$, -1 if edge $j$
is beginning with vertex $n_i$, and 0 otherwise. Define the {\it
Laplacian matrix} of $\cal{G}$ as $L{(\cal{G})}=BB^T.$ It follows
that $L{(\cal{G})}=\Delta-{\cal{A}},$ where ${\cal{A}}$ is the
adjacency matrix of undirected graph ${\cal{G}}$ and $\Delta$ is a
diagonal matrix whose $i$th diagonal element is the valence of
vertex $n_i$ in the graph. $L$ is always positive semi-definite.
Moreover, for a connected graph, $L$ has a single zero eigenvalue,
and the associated right eigenvector is ${\bf 1}_m$.

\end{document}